\documentclass{article}
\usepackage{epsfig}
\usepackage{amsfonts}
\usepackage{amsmath}

\newcommand{\field}[1]{\mathbb{#1}}
\newcommand{\R}{\field{R}}

\newcommand{\bmtx}[1]{\left[\begin{array}{#1}}
\newcommand{\emtx}{\end{array}\right]} 
\newcommand{\bsmtx}{\left[ \begin{smallmatrix}}
\newcommand{\esmtx}{\end{smallmatrix} \right]} 
\newcommand{\be}{\begin{equation}}
\newcommand{\ee}{\end{equation}}

\newcommand{\dv}{d}

\begin{document}
%-------------------------------------------------------
%       Title/Author Information
%-------------------------------------------------------
\title{SOSOPT: A Toolbox for Polynomial Optimization \\ Version 2.00}
\author{Pete Seiler \\ \texttt{seiler@aem.umn.edu} }
\date{\today}
\maketitle    

%-------------------------------------------------------
%                       Abstract
%-------------------------------------------------------
\begin{abstract}
  
  SOSOPT is a Matlab toolbox for formulating and solving
  Sum-of-Squares (SOS) polynomial optimizations. This document briefly
  describes the use and functionality of this toolbox.
  Section~\ref{sec:intro} introduces the problem formulations for SOS
  tests, SOS feasibility problems, SOS optimizations, and generalized
  SOS problems. Section~\ref{sec:using} reviews the SOSOPT toolbox for
  solving these optimizations.  This section includes information on
  toolbox installation, formulating constraints, solving SOS
  optimizations, and setting optimization options. Finally,
  Section~\ref{sec:sdp} briefly reviews the connections between SOS
  optimizations and semidefinite programs (SDPs).  It is the
  connection to SDPs that enables SOS optimizations to be solved in an
  efficient manner.

\end{abstract}

%-------------------------------------------------------
% Sum of Squares Optimization
%-------------------------------------------------------
\section{Sum of Squares Optimizations}
\label{sec:intro}

This section describes several optimizations that can be formulated
with sum-of-squares (SOS) polynomials
\cite{parrilo00,lasserre01,parrilo03}.  A multivariable polynomial is
a SOS if it can be expressed as a sum of squares of other polynomials.
In other words, a polynomial $p$ is SOS if there exists polynomials
$\{f_i\}_{i=1}^m$ such that $p = \sum_{i=1}^m f_i^2$. An SOS
polynomial is globally nonnegative because each squared term is
nonnegative.  This fact enables sufficient conditions for many
analysis problems to be posed as optimizations with polynomial SOS
constraints.  This includes many nonlinear analysis problems such as
computing regions of attraction, reachability sets, input-output
gains, and robustness with respect to uncertainty for nonlinear
polynomial systems
\cite{parrilo00,tibken00,hachicho02,wloszek03b,wloszek03a,parrilo03,sostools04,
  yalmip04,gatermann04,chesi04,tan04,
  wloszek05,papa05,tan06,tan06acc,tibken06,topcu07,topcu08a,topcu08b,tan08,
  topcu09,nasaworkshop09}.  The remainder of this section defines SOS
tests, SOS feasibility problems, SOS optimizations, and generalized SOS
optimizations.

Given a polynomial $p(x)$, a \underline{\bf sum-of-squares test}
is an analysis problem of the form:
\begin{align}
\label{eq:sostest}
\mbox{Is } p \mbox{ a SOS?}
\end{align}

A \underline{\bf sum-of-squares feasibility problem} is to construct
decision variables to ensure that certain polynomials are SOS.  More
specifically, an SOS feasibility problem is an optimization with
constraints on polynomials that are affine functions of the decision
variables:
\begin{align}
\label{eq:sosfeas}
& \mbox{ Find } \dv \in \R^r \mbox{ such that } \\
\nonumber
& a_k(x,\dv) \in \mbox{ SOS}, \ \ k=1,\ldots N_s \\
\nonumber
& b_j(x,\dv) = 0, \ \ j=1,\ldots N_e 
\end{align}
$\dv\in \R^r$ are decision variables.  The polynomials $\{ a_k \}$ and
$\{ b_j \}$ are given as part of the problem data and are affine in
$\dv$, i.e.  they are of the form:
\begin{align*}
a_k(x,\dv) & :=a_{k,0}(x) + a_{k,1}(x)\dv_1 + \dots + a_{k,n}(x)\dv_n \\
b_j(x,\dv) & :=b_{j,0}(x) + b_{j,1}(x)\dv_1 + \dots + b_{j,n}(x)\dv_n 
\end{align*}

A \underline{\bf sum-of-squares optimization} is a problem with a
linear cost and constraints on polynomials that are affine functions
of the decision variables:
\begin{align}
\label{eq:sosprog}
& \min_{\dv\in\R^r} c^T \dv\\
\nonumber
& \mbox{subject to: } \\
\nonumber
& a_k(x,\dv) \in \mbox{ SOS}, \ \ k=1,\ldots N_s \\
\nonumber
& b_j(x,\dv) = 0, \ \ j=1,\ldots N_e 
\end{align}
Again, $\dv\in \R^r$ denotes the decision variables and the
polynomials $\{ a_k \}$ and $\{ b_j \}$ are given polynomials that are
affine in $\dv$. SOS tests, feasibility problems, and optimizations
are all convex optimization problems.  These problems are solved by
exploiting the connections between SOS polynomials and positive
semidefinite matrices.  This is briefly reviewed in the
Section~\ref{sec:sdp}.

Finally, a \underline{\bf generalized sum-of-squares optimization} is
a problem of the form:
\begin{align}
\label{eq:gsosprog}
& \min_{\dv\in\R^r,t\in \R} t \\
\nonumber
& \mbox{subject to: } \\
\nonumber
& t b_k(x,\dv) - a_k(x,\dv) \in \mbox{ SOS}, \ \ k=1,\ldots N_g \\
\nonumber
& b_k(x,\dv) \in \mbox{ SOS}, \ \ k=1,\ldots N_g \\
\nonumber
& c_j(x,\dv) = 0, \ \ j=1,\ldots N_e 
\end{align}
$t\in\R$ and $\dv\in \R^r$ are decision variables.  The polynomials
$\{ a_k \}$, $\{ b_k \}$, and $\{ c_k \}$ are given data and are
affine in $\dv$.  The optimization cost is $t$ which is linear in the
decision variables.  The optimization involves standard SOS and
polynomial equality constraints. However, this is not an SOS
optimization because the constraints, $t b_k(x,\dv) - a_k(x,\dv$ is
SOS, are bilinear in the decision variables $t$ and $u$.  However, the
generalized SOS program is quasiconvex \cite{seiler10} and it can also
be solved efficiently as described in the next subsection.

%-------------------------------------------------------
% SOSOPT
%-------------------------------------------------------
\section{Using SOSOPT}
\label{sec:using}

This section describes the \texttt{sosopt} toolbox for solving SOS
optimizations.

%  Installation
\subsection{Installation}
\label{subsec:install}

The toolbox was tested with MATLAB versions R2009a and R2009b.  To
install the toolbox:
\begin{itemize}
\item Download the zip file and extract the contents to the directory
  where you want to install the toolbox.  
\item Add the \texttt{sosopt} directory to the Matlab path, e.g. using
  Matlab's \texttt{addpath} command.  
\end{itemize}

The \texttt{sosopt} toolbox requires the \texttt{multipoly} toolbox to
construct the polynomial constraints.  \texttt{multipoly} can be
obtained from \texttt{http://www.aem.umn.edu/$\sim$AerospaceControl/}.
\texttt{sosopt} also requires one of the following optimization codes
for solving semidefinite programs (SDPs): SeDuMi, SDPT3, CSDP,
DSDP, SDPAM, or SDPLR.  \texttt{sosopt} has been most extensively
tested on SeDuMi version 1.3 \cite{sedumi99,sturm01}.  The latest
version of SeDuMi can be obtained from
\texttt{http://sedumi.ie.lehigh.edu/}.

% Constraints
\subsection{Formulating Constraints}
\label{subsec:polyconstr}

Polynomial SOS and equality constraints are formulated using
\texttt{multipoly} toolbox objects.  The relational operators
\texttt{<=} and \texttt{>=} are overloaded to create SOS constraints.
If $p$ and $q$ are polynomials then \texttt{p>=q} and \texttt{p<=q}
denote the constraints $p-q \in$ SOS and $q-p \in$ SOS,
respectively.  The relational operator \texttt{==} is overloaded to
create a polynomial equality constraint.  If $p$ and $q$ are
polynomials then \texttt{p==q} denotes the constraint $p-q=0$. These
overloaded relational operators create a \texttt{polyconstr}
constraint object.  For example, the following code constructs the
constraints $6+d_1 x_1^2 - 5x_2^2 \in$ SOS and $d_1 x_1^2+d_2-
6x_1^2+4=0$.
\begin{verbatim}
>> pvar x1 x2 d1 d2
>> p = 6+d1*x1^2;
>> q = 5*x2^2;
>> p>=q
ans = 
  d1*x1^2 - 5*x2^2 + 6
  >= 0

>> class(ans)
ans =
polyconstr

>> p=d1*x1^2+d2;
>> q=6*x1^2+4
>> p==q
ans = 
  d1*x1^2 - 6*x1^2 + d2 - 4
  == 0
\end{verbatim}
The polynomial constraints are displayed in a standard form with all
terms moved to one side of the constraint.  The polynomials on the
left and right sides of the constraint are stored and can be accessed
with \texttt{.LeftSide} and \texttt{.RightSide}. The one-sided
constraint that is displayed can be accessed with \texttt{.OneSide}.
In addition, multiple polynomial constraints can be stacked into a
vector list of constraints using the standard Matlab vertical
concatenation with brackets and rows separated by a semicolon.
Finally, it is also possible to reference and assign into a list of
polynomial constraints using standard Matlab commands. These features
are shown below.
\begin{verbatim}
>> pvar x1 x2 d1 d2
>> constraint1 =  6+d1*x1^2 >= 5*x2^2;
>> constraint1.LeftSide
ans = 
  d1*x1^2 + 6

>> constraint1.RightSide
ans = 
  5*x2^2

>> constraint1.OneSide
ans = 
  d1*x1^2 - 5*x2^2 + 6

>> constraint2 =  d1*x1^2+d2 == 6*x1^2+4;
>> constraints = [constraint1; constraint2]
constraints = 
  polyconstr object with 2 constraints.

>> constraints(1)
ans = 
  d1*x1^2 - 5*x2^2 + 6
  >= 0

>> constraints(1).OneSide
ans = 
  d1*x1^2 - 5*x2^2 + 6

>> constraints(2)
ans = 
  d1*x1^2 - 6*x1^2 + d2 - 4
  == 0

>> constraints(2) = (d2==8);
>> constraints(2)
ans = 
  d2 - 8
  == 0

>> constraints.RelOp
ans = 
    '>='
    '=='

\end{verbatim}

% Solving
\subsection{Solving SOS Optimizations}
\label{subsec:solve}

The four SOS problems introduced in Section~\ref{sec:intro} can
be solved using the \texttt{sosopt} functions described below.
Documentation for each function can be obtained at the Matlab
prompt using the \texttt{help} Command.

\begin{enumerate}

\item \underline{\bf SOS test}: The function \texttt{issos} tests
if a polynomial $p$ is SOS.  The syntax is:
\begin{verbatim}
[feas,z,Q,f] = issos(p,opts)
\end{verbatim}

\texttt{p} is a \texttt{multipoly} polynomial object.  \texttt{feas} is
equal to 1 if the polynomial is SOS and 0 otherwise.  If
\texttt{feas}=1 then $f$ is a vector of polynomials that provide
the SOS decomposition of $p$, i.e. $p = \sum_i f_i^2$. 
$z$ is a vector of monomials and and $Q$ is a positive semidefinite
matrix such that $p=z^TQz$. $z$ and $Q$ are a Gram matrix decomposition
for $p$. This is described in more detail in Section~\ref{sec:sdp}.
The \texttt{opts} input is an \texttt{sosoptions} object.  Refer to
Section~\ref{subsec:opts} for more details on these options.

\vspace{0.1in}

\item \underline{\bf SOS feasibility}: The function \texttt{sosopt}
solves SOS feasibility problems. The syntax is:
\begin{verbatim}
[info,dopt,sossol] = sosopt(pconstr,x,opts);
\end{verbatim}

\texttt{pconstr} is an $N_p \times 1$ vector of polynomial SOS and
equality constraints constructed as described in
Section~\ref{subsec:polyconstr}.  \texttt{x} is a vector list of
polynomial variables. The variables listed in \texttt{x} are the
independent polynomial variables in the constraints. All other
variables that exist in the polynomial constraints are assumed to be
decision variables. The polynomial constraints must be affine
functions of these decision variables.  The \texttt{opts} input is an
\texttt{sosoptions} object (See Section~\ref{subsec:opts}).

The \texttt{info} output is a structure that contains a variety of
information about the construction of the SOS optimization problem.
The main data in this structure is the \texttt{feas} field.  This
field is equal to 1 if the problem is feasible and 0 otherwise.

The \texttt{dopt} output is a polynomial array of the optimal decision
variables.  The first column of \texttt{dopt} contains the decision
variables and the second column contains the optimal values.  The
polynomial \texttt{subs} command can be used to replace the decision
variables in any polynomial with their optimal values, e.g
\texttt{subs( pconstr(1).LeftSide, dopt) } substitutes the optimal
decision variables into the left side of the first constraint.
\texttt{dopt} is returned as empty if the optimization is infeasible.

\texttt{sossol} is an $N_p \times 1$ structure array with fields
\texttt{p}, \texttt{z}, and \texttt{Q}.  \texttt{sossol(i).p} is
\texttt{pconstr(i)} evaluated at the optimal decision variables.  If
\texttt{pconstr(i)} is an SOS constraint then \texttt{sossol(i).z} and
\texttt{sossol(i).Q} are the vector of monomials and positive
semidefinite matrix for the Gram matrix decomposition of
\texttt{sossol(i).p}, i.e. $p=z^TQz$.  This Gram matrix decomposition
is described in more detail in Section~\ref{sec:sdp}.  If
\texttt{pconstr(i)} is a polynomial equality constraint then these two
fields are returned as empty.  \texttt{sossol} is empty if the
optimization is infeasible.

\vspace{0.1in}

\item \underline{\bf SOS optimization}: The function \texttt{sosopt}
also solves SOS optimization problems. The syntax is:
\begin{verbatim}
[info,dopt,sossol] = sosopt(pconstr,x,obj,opts);
\end{verbatim}

\texttt{obj} is a polynomial that specifies the objective function.
This must be be an affine function of the decision variables and it
cannot depend on the polynomial variables. In other words, \texttt{obj} 
must have the form $c_0 + \sum_i c_i d_i$ where $c_i$ are real numbers
and $d_i$ are decision variables. The remaining inputs and outputs
are the same as described for SOS feasibility problems.  The \texttt{info}
output has one additional field \texttt{obj} that specifies the
minimal value of the objective function.  This field is the same
as \texttt{subs(obj,dopt)}.  \texttt{obj} is set to \texttt{+inf}
if the problem is infeasible.

\vspace{0.1in}
  
\item \underline{\bf Generalized SOS optimization}: The function
  \texttt{gsosopt} solves generalized SOS optimization problems. The
  syntax is:
\begin{verbatim}
[info,dopt,sossol] = gsosopt(pconstr,x,t,opts)
\end{verbatim}

\texttt{pconstr} is again an $N_p \times 1$ vector of polynomial SOS and
equality constraints constructed as described in
Section~\ref{subsec:polyconstr}.  \texttt{x} is a vector list of
polynomial variables. The variables listed in \texttt{x} are the
independent polynomial variables in the constraints. All other
variables that exist in the polynomial constraints are assumed to be
decision variables.  The objective function
is specified by the third argument \texttt{t}. This objective must be a 
single polynomial variable and it must be one of the decision variables.
The constraints must have the special structure specified in the 
Generalized SOS problem formulation.  Let (\texttt{d},\texttt{t}) 
denote the complete list of decision variables.
The constraints are allowed to have bilinear terms involving products
of \texttt{t} and \texttt{d}.  However, they must be linear in \texttt{d}
for fixed \texttt{t} and linear in \texttt{t} for fixed \texttt{d}.
The \texttt{opts} input is an \texttt{gsosoptions} object (See
Section~\ref{subsec:opts}).

The outputs are the same as described for SOS feasibility and
optimization problems.  The only difference is that the 
\texttt{info} output does not have an \texttt{obj} field. \texttt{gsosopt}
uses a bisection to solve the generalized SOS problem. It 
computes lower and upper bounds on the optimal cost such that the
bounds are within a specified stopping tolerance.  These bounds are
returned in the \texttt{tbnds} field.  This is a $1 \times 2$ vector
[$t_{lb}$, $t_{ub}$] giving the lower bound $t_{lb}$ and upper bound
$t_{ub}$ on the minimum value of \texttt{t}. \texttt{tbnds} is empty 
if the optimization is infeasible.

\end{enumerate}

% Solving
\subsection{Constructing Polynomial Decision Variables}
\label{subsec:decvar}

The \texttt{sosopt} and \texttt{multipoly} toolboxes contain several
functions to quickly and easily construct polynomials whose
coefficients are decision variables.  The \texttt{mpvar} and
\texttt{monomials} functions in the \texttt{multipoly} toolbox can be
used to construct a matrix of polynomial variables and a vector list
of monomials, respectively. Examples are shown below:
\begin{verbatim}
>> P = mpvar('p',[4 2])
P = 
  [ p_1_1, p_1_2]
  [ p_2_1, p_2_2]
  [ p_3_1, p_3_2]
  [ p_4_1, p_4_2]

>> pvar x1 x2
>> w = monomials([x1;x2],0:2)
w = 
  [     1]
  [    x1]
  [    x2]
  [  x1^2]
  [ x1*x2]
  [  x2^2]

\end{verbatim}
The first argument of \texttt{mpvar} specifies the prefix for the
variable names in the matrix and the the second argument specifies the
matrix size.  The first argument of \texttt{monomials} specifies the variables
used to construct the monomials vector.  The second argument specifies
the degrees of monomials to include in the monomials vector. In the
example above, the vector \texttt{w} returned by \texttt{monomials}
contains all monomials in variables \texttt{x1} and \texttt{x2} of
degrees 0,1, and 2.

These two functions can be used to quickly construct a polynomial $p$
that is a linear combination of monomials in $x$ with
coefficients specified by decision variables $d$.
\begin{verbatim}
>> pvar x1 x2
>> w = monomials([x1;x2],0:2);
>> d = mpvar('d',[length(w),1]);
>> [w, d]
ans = 
  [     1, d_1]
  [    x1, d_2]
  [    x2, d_3]
  [  x1^2, d_4]
  [ x1*x2, d_5]
  [  x2^2, d_6]

>> p = d'*w
p = 
  d_4*x1^2 + d_5*x1*x2 + d_6*x2^2 + d_2*x1 + d_3*x2 + d_1

\end{verbatim}
This example constructs a quadratic function in variables $(x_1,x_2)$
with coefficients given by the entries of $d$. $p$ could alternatively
be interpreted as a cubic polynomial in variables $(x,d)$.

The \texttt{polydecvar} function can be used to construct
polynomials of this form in one command:
\begin{verbatim}
>> p = polydecvar('d',w)
p = 
  d_4*x1^2 + d_5*x1*x2 + d_6*x2^2 + d_2*x1 + d_3*x2 + d_1
\end{verbatim}
The first argument of \texttt{polydecvar} specifies the prefix for the
coefficient names and the second argument specifies the monomials to
use in constructing the polynomial. The output of \texttt{polydecvar}
is a polynomial in the form: \texttt{p=d'*w} where \texttt{d} is a
coefficient vector generated by \texttt{mpvar}.  This is called the
\underline{vector} form because the coefficients are specified in the
vector \texttt{d}.

The Gram matrix provides an alternative formulation for specifying
polynomial decision variables.  In particular, one can specify a
polynomial as $p(x,D)=z(x)^T D z(x)$ where $z(x)$ is a vector of
monomials and $D$ is a symmetric matrix of decision variables.  A
quadratic function in variables $(x_1,x_2)$ with coefficient matrix
$D$ is constructed as follows:
\begin{verbatim}
>> pvar x1 x2
>> z = monomials([x1;x2],0:1);
>> D = mpvar('d',[length(z) length(z)],'s')
D = 
  [ d_1_1, d_1_2, d_1_3]
  [ d_1_2, d_2_2, d_2_3]
  [ d_1_3, d_2_3, d_3_3]

>> s = z'*D*z
s = 
  d_2_2*x1^2 + 2*d_2_3*x1*x2 + d_3_3*x2^2 + 2*d_1_2*x1 + 2*d_1_3*x2 + d_1_1

\end{verbatim}
The \texttt{'s'} option specifies that \texttt{mpvar} should return
a symmetric matrix.  This construction can be equivalently performed
using the \texttt{sosdecvar} command:
\begin{verbatim}
>> [s,D] = sosdecvar('d',z)
s = 
  d_2_2*x1^2 + 2*d_2_3*x1*x2 + d_3_3*x2^2 + 2*d_1_2*x1 + 2*d_1_3*x2 + d_1_1
D = 
  [ d_1_1, d_1_2, d_1_3]
  [ d_1_2, d_2_2, d_2_3]
  [ d_1_3, d_2_3, d_3_3]
\end{verbatim}
This is called the \underline{matrix} form because the coefficients
are specified in the symmetric matrix \texttt{D}.

In the examples above, the vector and matrix forms both use six
independent coefficients to specify a quadratic polynomial in
$(x_1,x_2)$. In general, the matrix form uses many more variables than
the vector form to represent the coefficients of a polynomial.  Thus
the vector form will typically lead to more efficient problem
formulations.  The only case in which \texttt{sosdecvar} leads to more
efficient implementations is when the resulting polynomial is directly
constrained to be SOS.  Specifically, the \texttt{sosdecvar} command
should be used to construct polynomials that will be directly added to
the list of SOS constraints, as in the example below:
\begin{verbatim}
>> [s,D] = sosdecvar('d',z);
>> pconstr(i) = s>=0;
\end{verbatim}
\textbf{NOTE: Creating a polynomial variable \texttt{s} using
  the \texttt{sosdecvar} command will not cause \texttt{sosopt} or
  \texttt{gsosopt} to constrain the polynomial to be SOS.  The constraint
  \texttt{s>=0} must be added to the list of constraints to enforce 
  \texttt{s} to be SOS.}

\subsection{Demos}
\label{subsec:demo}

\texttt{sosopt} includes several demo files that illustrate the use of
the toolbox.  These demo files can be found in the \texttt{Demos}
subfolder. A brief description of the existing demo files is given
below.

\begin{enumerate}
  
\item \underline{\bf SOS test}: \texttt{issosdemo1} demonstrates the
  use of the \texttt{issos} function for testing if a polynomial $p$ 
is a sum of squares.  This example uses \texttt{issos} to construct an
SOS decomposition for a degree four polynomial in two variables. The
example polynomial is taken from Section 3.1 of the SOSTOOLs
documentation \cite{sostools04}. \texttt{sosoptdemo1} solves the
same SOS test using the \texttt{sosopt} function.

\item \underline{\bf SOS feasibility}: There are three demo files that
  solve SOS feasibility problems: \texttt{sosoptdemo2},
  \texttt{sosoptdemo4}, and \texttt{sosoptdemo5}. These examples are
  taken from Sections 3.2, 3.4, and 3.5 of the SOSTOOLs documentation
  \cite{sostools04}, respectively. Demo 2 solves for a global Lyapunov
  function of a rational, nonlinear system. Demo 4 verifies the
  copositivity of a matrix. Demo5 computes an upper bound for
  a structured singular value problem.
  
\item \underline{\bf SOS optimization}: There are three demo files
  that solve SOS optimization problems: \texttt{sosoptdemo3},
  \texttt{sosoptdemoLP}, and \texttt{sosoptdemoEQ}. Demo 3 is taken
  from Section 3.3 of the SOSTOOLs documentation \cite{sostools04}.
  This demo uses SOSOPT to compute a lower bound on the global minimum
  of the Goldstein-Price function.  The EQ demo provides a simple
  example with polynomial equality constraints in addition to SOS
  constraints.  Finally, the LP demo shows that linear programming
  constraints can be formulated using \texttt{sosopt}.
  
\item \underline{\bf Generalized SOS optimization}: There are two demo
  files that solve generalized SOS optimization problems:
  \texttt{gsosoptdemo1} and \texttt{pcontaindemo1}. \texttt{gsosoptdemo1}
  \texttt{gsosopt} to compute an estimate of the region of 
  attraction for the van der Pol oscillator using the Lyapunov 
  function obtained via linearization. \texttt{pcontaindemo1}
  solves for the radius of the largest circle that lies within the
  contour of a 6th degree polynomial.  This is computed using the
  specialized function \texttt{pcontain} for verifying set containments.
  The set containment problem is a specific type of generalized
  SOS optimization.

\end{enumerate}

% Solving
\subsection{Options}
\label{subsec:opts}

The \texttt{sosoptions} command will create a default options structure for
the \texttt{issos} and \texttt{sosopt} functions. 
The \texttt{sosoptions} command will return an object with 
the fields:
\begin{itemize}
\item \texttt{solver}: Optimization solver to be used. The choices
  are: 'sedumi', 'sdpam', 'dsdp', 'sdpt3', 'csdp', or 'sdplr'.
  The default solver is 'sedumi'.
\item \texttt{form}: Formulation for the optimization. The choices are
  'image' or 'kernel'. These forms are described in
  Section~\ref{sec:sdp}.  The default is 'image'.
\item \texttt{simplify}: SOS simplification procedure to remove
  monomials that are not needed in the Gram matrix form. This reduces
  the size of the related semidefinite programming problem and hence
  also reduces the computational time. The choices are 'on' or 'off'
  and the default is 'on'.
\item \texttt{scaling}: Scaling of SOS constraints. This scales each
  constraint by the Euclidean norm (2-norm) of the one-sided
  polynomial coefficient vector.  The choices are 'on' or 'off' and
  the default is 'off'.
\item \texttt{checkfeas}: Check feasibility of solution.  The choices
  are 'off', 'fast', 'full', and 'both'.  The default is 'fast'.
  'fast' checks feasibility information returned by the solver. 'full'
  checks the validity of the Gram matrix decomposition in the output
  sossol. The 'full' check is more computationally costly. 'both' does
  both feasibility checks.
\item \texttt{feastol}: Feasibility tolerance used in the 'full'
  feasibility check.  This should be a positive, scalar, double. The
  default is \texttt{1e-6}. 
\item \texttt{solveropts}: Structure with options passed directly to
  the optimization solver.  The default is empty.  The solver display
  is turned off with this default.
\end{itemize}

\vspace{0.1in}

The \texttt{gsosoptions} command will create a default options
structure for the \texttt{gsosopt} function.  The \texttt{gsosoptions}
command will return an object with all fields contained in an
\texttt{sosoptions} structure.  In addition it will contain the
fields:
\begin{itemize}
\item \texttt{minobj}: Minimum value of objective for bisection. 
This should be a scalar double. The default is \texttt{-1e3}. 
\item \texttt{maxobj}: Maximum value of objective for bisection. 
This should be a scalar double. Moreover, \texttt{maxobj} should
be $\ge$ \texttt{minobj}. The default is \texttt{1e3}. 
\item \texttt{absbistol}: Absolute bisection stopping tolerance This
  should be a positive, scalar, double. The default is \texttt{1e-3}.
  The bisection terminates if $t_{ub}-t_{lb} \le$ \texttt{absbistol}.
\item \texttt{relbistol}: Relative bisection stopping tolerance This
  should be a positive, scalar, double. The default is \texttt{1e-3}.
  The bisection terminates if $t_{ub}-t_{lb} \le$ \texttt{relbistol}
  $\times$ $t_{lb}$.
\item \texttt{display}: Display bisection iteration information.  The
  choices are 'on' or 'off' and the default is 'off'.  If
  \texttt{display = 'on'} then \texttt{gsosoptions} displays, for each
  iteration, the attempted value of $t$, feasibility result and the
  current upper and lower bounds on the optimal value of $t$.  The
  display information generated by the optimization solver is not
  affected by this option.
\end{itemize}

%-------------------------------------------------------
% Connections to SDPs
%-------------------------------------------------------
\section{Connections to SDPs}
\label{sec:sdp}

Given a polynomial $p$, an SOS test is to determine if $p$ is a SOS.
To solve this problem, the polynomial if expressed in the form $p =
z^TQz$ where $z$ is a vector of monomials and $Q$ is a symmetric
``Gram'' matrix \footnote{ A monomial is a term of the form
  $x^\alpha\doteq x_1^{\alpha_1} x_2^{\alpha_2} \cdots x_n^{\alpha_n}$
  where the $\alpha_i$ are non-negative integers.}.  The Gram matrix
is not unique and a known result is that $p$ is a SOS if and only if
there exists $Q=Q^T \succeq 0$ such that $p = z^T Q z$
\cite{choi95,powers98}. Equating the coefficients of $p$ and $z^TQz$
leads to linear equality constraints on the entries of $Q$.  There
exists a matrix $A$ and vector $b$ such that these equality
constraints can be represented as $Aq = b$ where $q:=vec(Q)$ denotes
the vector obtained by vertically stacking the columns of $Q$. Thus
the SOS test can be converted to a problem of the form: 
\begin{align}
\label{eq:image}
  \mbox{Given a matrix $A$ and vector $b$, find $Q \succeq 0$ such
    that $Aq = b$.}
\end{align}
This is a semidefinite programming (SDP) problem
\cite{befb94,vandenberghe96}. In general, there are fewer equality
constraints than independent entries of $Q$, i.e.  $A$ has fewer rows
than columns.  One can compute a particular solution $Q_0$ such that
$p=z^TQ_0z$ and a basis of homogeneous solutions $\{N_i\}$ such that
$z^TN_iz = 0$ for each $i$ where $0$ is the zero polynomial.
The matrix $A$ has special structure that can be exploited to
efficiently compute these matrices. Thus every matrix $Q$ satisfying
$p=z^TQz$ can be expressed in the form $Q_0 + \sum_i \lambda_i N_i
\succeq 0$ where $\lambda_i \in \R$. This enables the
SOS test to be converted into the alternative formulation:
\begin{align}
\label{eq:kernel}
  \mbox{Given matrices $Q_0$ and $\{N_i\}$ find a vector $\lambda$
    such that $Q_0 + \sum_i \lambda_i N_i \succeq 0$.}
\end{align}
This problem has a single linear matrix inequality (LMI) and is also a
semidefinite programming problem.  The SDPs in Equation \ref{eq:image}
and Equation \ref{eq:kernel} are dual optimization problems
\cite{vandenberghe96}.  There exist many freely available codes to
solve these types of problem, e.g.  SeDuMi \cite{sedumi99,sturm01}.
In the SeDuMi formulation, Equation~\ref{eq:image} is called the
primal or image problem and Equation~\ref{eq:kernel} is the dual or
kernel problem.

The constraints in SOS feasibility and optimization problems are
similarly converted to semidefinite matrix constraints. For example,
$a_k(x,d)$ is SOS if and only if there exists $Q \succeq 0$ such
that
\begin{align}
a_{k,0}(x) + a_{k,1}(x)\dv_1 + \dots + a_{k,n}(x)\dv_n  = z(x)^T Q z(x)
\end{align}
Equating the coefficients leads to linear equality constraints on the
decision variables $\dv$ and the entries of $Q$.  There exist matrices
$A_\dv$, $A_q$ and a vector $b$ such that these equality constraints
can be represented as $A_\dv \dv + A_q q = b$ where $q:=vec(Q)$. Thus
$a_k(x,d)$ is SOS if and only if there exists $Q \succeq 0$ such that
$A_\dv \dv + A_q q = b$.  Each SOS constraint can be replaced in this
way by a positive semidefinite matrix subject to equality constraints
on its entries and on the decision variables.  The polynomial equality
constraints are equivalently represented by equality constraints on
the decision variables.  Performing this replacement for each
constraint in an SOS feasibility or optimization problem leads to an
optimization with equality and semidefinite matrix constraints.  This
is an SDP in SeDuMi primal/image form.  An SDP in SeDuMi dual/kernel
is obtained by replacing the positive semidefinite matrix variables
$Q$ that that arise from each SOS constraint with linear combinations
of a particular solution $Q_0$ and homogeneous solutions $\{N_i\}$.
This is similar to the steps described above for the SOS test and full
details can be found in \cite{nasaworkshop09}.

Finally, the generalized SOS optimization has SOS constraints that are
bilinear in decision variables $t$ and $d$.  A consequence of this
bilinearity is that the SOS constraints cannot be replaced with linear
equality constraints on the decision variables. However, the
generalized SOS program is quasiconvex \cite{seiler10} and it can be
efficiently solved.  In particular, for fixed values of $t$ the
constraints are linear in the remaining decision variables $\dv$.  An
SOS feasibility problem can be solved to determine if the constraints
are feasible for fixed $t$.  Bisection can be used to find the minimum
value of $t$, to within a specified tolerance, for which the
constraints are feasible.  In principle this problem can also be
converted to a generalized eigenvalue problem \cite{boyd93} (subject
to some additional technical assumptions) but the theory and available
software for generalized eigenvalue problems are not as well-developed
as for SDPs.

\texttt{sosopt} converts the SOS optimizations into SDPs in either
primal/image or dual/kernel form.  The form can be specified with the
\texttt{form} option in the \texttt{sosoptions} object.  Interested
users can see the lower level functions \texttt{gramconstraint} and
\texttt{gramsol} for implementation details on this conversion.
\texttt{sosopt} then solves the SDP using one of the freely available
solvers that have been interfaced to the toolbox.  The \texttt{solver}
option is used to specify the solver.  Finally, \texttt{sosopt}
converts the SDP solution back to polynomial form.  Specifically, the
optimal SOS decision variables and the Gram matrix decompositions are
constructed from the SDP solution.  \texttt{sosopt} also checks the
feasibility of the returned solution.  The \texttt{checkfeas} option
specifies the feasibility check performed by \texttt{sosopt}.  The
\texttt{fast} option simply checks the feasibility information
returned by the SDP solver. The \texttt{full} option verifies the Gram
matrix decomposition for each SOS constraint.  In particular, it
checks that the Gram matrix is positive semidefinite and it checks
that $p=z^TQz$ within some tolerance.  The \texttt{full} feasibility
check also verifies that each SOS equality constraint is satisfied
within a specified tolerance.

%-------------------------------------------------------
% Acknowledgments
%-------------------------------------------------------
\section{Acknowledgments}

This research was partially supported under the NASA Langley NRA
contract NNH077ZEA001N entitled ``Analytical Validation Tools for Safety
Critical Systems'' and the NASA Langley NNX08AC65A contract 
entitled 'Fault Diagnosis, Prognosis and Reliable Flight 
Envelope Assessment.'' The technical contract monitors are Dr. Christine
Belcastro and Dr. Suresh Joshi, respectively.

%-------------------------------------------------------
%                REFERENCES (IN .BIB FILE)             
%-------------------------------------------------------
\bibliography{sosoptdoc}

\end{document}